\input amstex.tex
 \documentstyle{amsppt} 
 \magnification\magstep1

\hoffset=0cm
\voffset=1cm
\advance\vsize -3cm

 \document
 \topmatter
 \title  From gauge anomalies to gerbes and gerbal actions  \endtitle
 \author Jouko Mickelsson \endauthor
 \affil Department of Mathematics and Statistics, University of Helsinki, and
 Department of Theoretical Physics, Royal Institute of Technology, Stockholm 
\endaffil
 
 \endtopmatter

\define\Aut{\text{Aut}}
\define\Out{\text{Out}}

 ABSTRACT  The purpose of this contribution is to point out connections between 
 recent ideas about gerbes and gerbal actions (as higher categorical extension
 of representation theory)  and old discussion in quantum field theory on commutator
 anomalies, gauge group extensions, and 3-cocycles. The unifying concept is the
 classical obstruction theory for group extensions as explained in the
 reference [ML].
[Based on talk given at the meeting ``Motives, Quantum Field Theory,
 and Pseudodifferential operators'', Boston University, June 2-13, 2008]

 \vskip 0.4in
 
 1. INTRODUCTION 
 
 \vskip 0.3in
 
 It was first realized through perturbative analysis of gauge theories that the gauge 
 symmetry is broken in presence of chiral fermions, [ABJ]. Later, it was found that 
 this phenomen is related to index theory of (families) of Dirac operators. In
 particular, the effective action functional, defined as a regularized determinant of
 the Dirac operator, is not always gauge invariant and the lack of invariance can
 be formulated as the curvature of a complex line bundle, the determinant 
 line bundle, over the moduli space of gauge connections, [AS]. 
 
 In the hamiltonian formulation of gauge theory the symmetry breaking manifests
 itself as a modification of the commutation relations of the Lie algebra
 of infinitesimal gauge transformations. The gauge algebra (in case of trivial vector
 bundles) is the Lie algebra of functions $M\frak{g}$ from the physical space $M$ to
 a finite-dimensional Lie algebra $\frak g.$  The commutation relations of the modified 
 algebra can be written as 
 $$[(X, a), (Y,b)] = ([X,Y], c(A;X,Y))$$ 
 where $[X,Y]$ is the point-wise commutator in $M\frak g$ and $a,b$  are complex valued 
 functions of the gauge potential $A.$  $c$ is a Lie algebra 2-cocycle determining an extension
 of $M\frak g.$ In the case when $M$ is the unit circle $S^1$ it turns out that $c$ is independent
 of $A$ and we have a central extension defining (when $\frak g$ is simple) an affine
 Kac-Moody algebra. 
 
 When dim$M >1$ the cocycle $c$ depends explicitly on $A.$ It is still an open question 
 whether this algebra has interesting faithful Hilbert space representations, analogous to
 the highest weight representations of affine Lie algebras extensively used in string theory and constructions 
 of quantum field theory models in 1+1 space-time dimensions.  What is know at present that there are 
 natural  unitary \it Hilbert bundle \rm actions of the extended gauge Lie algebras and groups. These 
 come from quantizing chiral fermions in background gauge fields.  For each gauge connection $A$ there
 is a fermionic Fock space $\Cal F_A$ where the quantized Dirac hamiltonian $\hat D_A$ acts. This family 
 of (essentially positive) Dirac operators transforms equivariantly with respect to the action of the extension
 of the group $MG.$  The gauge transformations are defined as projective unitary operators between the fibers 
 $\Cal F_A$ and $\Cal F_{A^g}$ of the Fock bundle, corresponding to a true action of the gauge group 
 extension. As a consequence, the Fock bundle is defined only as a projective bundle over the 
 moduli space $\Cal A/ MG.$  Actually, in order that the moduli space is a smooth manifold, one has 
 to restrict $MG$ to the \it based gauge transformations \rm which are functions on $M$ taking the value 
 $e\in G$ at a fixed base point $x_0\in M.$ 
 
 A projective bundle is completely determined, up to equivalence, by the Dixmier-Douady class which 
 is an element of $H^3(\Cal A/ MG, \Bbb Z).$ This is the origin of gerbes in quantum field theory,
 [CMM]. Topologically a gerbe on a space $X$ is just an equivalence class of $PU(H)$ bundles
 over $X.$   In terms of \v{C}ech cohomology subordinate to a good cover $\{U_{\alpha}\}$ of $X,$ 
 the gerbe is given as a $\Bbb C^{\times}$ valued cocycle $\{f_{\alpha\beta\gamma}\},$ 
 $$ f_{\alpha\beta\gamma} f_{\alpha\beta\delta}^{-1} f_{\alpha\gamma\delta} f_{\beta\gamma\delta}^{-1} = 1$$
 on intersections $U_{\alpha}\cap U_{\beta}\cap U_{\gamma} \cap U_{\delta}.$ This cocycle arises from 
 the lifting problem: A  $PU(H)$ bundle is given in terms of transition functions $g_{\alpha\beta}$ with
 values in $PU(H).$ After lifting these to $U(H)$ one gets a family of functions $\hat g_{\alpha\beta}$ which satisfy the 1-cocycle 
 condition up to a phase,
 $$ \hat g_{\alpha\beta} \hat g_{\beta\gamma} \hat g_{\gamma\alpha} = f_{\alpha\beta\gamma} \bold{1}.$$
 
 The notion of gerbal action was introduced in the recent paper [FZ]. This is to be viewed 
 as the next  level after projective actions related to central extensions of groups and is given in
 terms of third group cohomology. In fact, the appearance of third cohomology in this context is 
 not new and is related to group extensions as explained in [ML]. In the simples form, the problem 
 is the following. Let  $F$ be an extension of $G$ by the group $N,$
 $$1\to N \to F \to G \to 1$$ 
 an exact sequence of groups. Suppose that $1\to a \to \hat N \to N \to 1$ is a central extension by
 the abelian group $a.$ Then one can ask whether the extension $F$ of $G$ by $N$ can be prolonged to
 an extension of $G$ by the group $\hat N.$  An obstruction to this is an element in the group cohomology
 $H^3(G, a)$ with coefficients in $a.$  In case of Lie groups, there is a corresponding  Lie algebra cocycle 
 representing a class in   $H^3(\frak g, a)$ We shall demonstrate this in detail for an example arising from
 quantization of gauge theory. It is closely related to the idea  in
 [Ca], further elaborated in [CGRS],   which in turn was a response 
 to a discussion in the 80's on breaking of Jacobi identity for the field algebra in Yang-Mills theory
 [GJJ].
 
 The paper is organized as follows. In Section 2 we explain the gauge group extensions arising from action
 on bundles of fermionic Fock spaces over background gauge fields and the corresponding Lie algebra cocycles .
 Section 3 consists of a general discussion how gerbal action arises from the group of outer automorphisms 
 of an associative algebra, how this leads to a 3-cocycle on the symmetry group, and finally we give
 an example coming from Yang-Mills theory in 1+1 space-time dimensions. Section 4 contains a generalization
 to Yang-Mills theory in higher space-time dimensions. Finally in Section 5 we explain an application to 
 twisted K-theory on moduli space of gauge connections.

 \vskip 0.3in
 
 2.  BUNDLES OF FOCK SPACES OVER GAUGE CONNECTIONS
 
 \vskip 0.3in
 
 A basic problem in quantum field theory in higher than two space-time dimensions
 is that the representations of canonical anticommutation relations algebra (CAR) 
 are not equivalent in different background gauge fields, and this leads to various
 divergencies in perturbation theory. However, in the case of the linear problem of 
 quantizing fermions in a background gauge field one can construct the hamiltonian and 
 the Hilbert space in a nonperturbative way. One can actually avoid the divergencies by
 taking systematically into account the need of dealing with a family of nonequivalent
 CAR algebra representations. 
 
 The method introduced in [Mi93] and generalized in [LM] is based on the observation 
 that for each gauge connection $A$ in the family $\Cal A$ of all gauge connections 
 on a vector bundle $E$ over a compact spin manifold, one can choose a unitary operator
 $T_A$ in the Hilbert space $H$ of $L^2$ sections in the tensor product of the spin bundle and
 the vector bundle $E$ such that the Dirac hamiltonian $D_A$ is conjugated to 
 $\tilde D_A=  T_A D_A T_A^{-1}$ such that the equivalent hamiltonian $\tilde D_A$ can 
 be quantized in the "free" Fock space, the Fock space for a fixed background connection
 $A_0.$ In the case of a trivial bundle $E$ one can take as $A_0$ the globally defined 
 gauge connection represented by 1-form equal to zero.
 
 The action of the group $\Cal G$ of smooth gauge transformations $A\mapsto A^g=g^{-1} A g
 +g^{-1}dg$ on the family $\tilde D_A$ is then given by
 $$\tilde D_A \mapsto \omega(A;g)^{-1} \tilde D_A \omega(A;g)$$
 corresponding to $D_A \mapsto g^{-1} D_A g$
 where $\omega(A;g) = T_{ A} g T_{A^ g}^{-1}$ satisfies the 1-cocycle relation 
 $$\omega(A;gg') = \omega(A;g) \omega(A^g;g').$$ 
 Furthermore, the cocycle satisfies the condition $[\epsilon, \omega(A;g)]$ is Hilbert-Schmidt.
 Here $\epsilon$ is the sign $D_{A_0}/|D_{A_0}| $ of the free Dirac operator. 
 This means that the operators $\omega(A;g)$ belong to the restricted unitary group
 $U_{res}(H_+\oplus H_-)$ where $H=H_+\oplus H_-$ is the polarization of  $H$
 with respect to the sign operator $\epsilon,$ [PS]. 
 
 The quantization of the operator $\tilde D_A$ is obtained in a fermionic Fock space $\Cal F$
 which carries an irreducible representation of the CAR algebra $\Cal B$ which is a completion
 of the algebra defined by generators and relations according to
 $$ a^*(u) a(v) + a(v)a^*(u)  = 2 <v,u>$$
 and all other anticommutators equal to zero. Here $u,v\in H$ and $<\cdot, \cdot>$ is the Hilbert
 space inner product (antilinear in the first argument). The representation is fixed (up to equivalence) by the requirement that
 there exists a vacuum vector $|0> \in \Cal F$ such that
 $$ a(u)|0> =0 = a^*(v)|0>, \text{ for } u\in H_+, v\in H_-.$$ 
 
 The group $U_{res}(H)$ has a central extension by $S^1$ such that the Lie algebra central extension
 is given by the 2-cocycle $c(X,Y) = \frac14 \text{tr\,}\epsilon[\epsilon,X] [\epsilon,Y].$ The central extension
 $\hat U_{res}$ has a unitary representation $g\mapsto \hat g$ in $\Cal F$ fixed by the requirement
 $$ \hat g a^*(u) g^{-1} = a^*(gu)$$ 
 for all $u\in H.$  
 
 If we choose a lift $\hat \omega(A;g)$ of  the element $\omega(A,g)$ to unitaries in the Fock space 
 $\Cal F$ we can write
 $$\hat \omega(A;gg') =\Phi(A; g,g') \hat\omega(A;g)\hat\omega(A^g;g')$$
 where $\Phi$ takes values in $S^1.$ It is a 2-cocycle by construction,
 $$\Phi(A;g,g')\Phi(A;gg',g'') = \Phi(A;g,g'g'') \Phi(A^g;g',g''),$$
 which is smooth in an open neighborhood of the neutral element in $\Cal G.$ 
 This just reflects the associativity in the group multiplication in the central extension $\hat U_{res}.$ 
 
 Taking the second derivative
 $$\frac{d^2}{dt ds}|_{t=s=0} \Phi(A; e^{tX}, e^{sY}) = \frac12 c(A;X,Y)$$ 
 gives a 2-cocycle  $c$ for the Lie algebra of $\Cal G$ with coefficients in the ring of complex functions of the 
 variable $A.$ 
 
 The cocycle depends on the lift $\omega\mapsto \hat \omega$ but two lifts are related by a multiplication by a circle
 valued function $\psi(A;g)$ and the corresponding 2-cocycles are related by a coboundary,
 $$\Phi'(A;g,g') =  \Phi(A,g,g') \psi(A;gg')\psi(A;g)^{-1} \psi(A^g;g')^{-1}.$$
 The Lie algebra cocycle $c$ satisfies
 $$ c(A;X, [Y,Z]) + \Cal L_X c(A;Y,Z) + \text{cyclic permutations} = 0,$$
 where $\Cal L_X$ is the Lie derivative acting on functions $f(A)$ through infinitesimal
 gauge transformation, $(\Cal L_X f)(A)= Df(A) \cdot ([A,X] +dX).$ 
 
 Explicite expressions for the cocycle $c$ have been computed in the literature, for example if the 
 physical space is a circle we get the central extension of a loop algebra (affine Kac-Moody algebra),
 $$c(A;X,Y) = \frac{1}{2\pi} \int_{S^1} \text{tr} \,XdY,\tag2.1$$
 where the trace is evaluated in a finite dimensional representation of $G.$ In this case $c$ does not
 depend on $A$ and the abelian extension reduces to a central extension. This reflects the fact that elements of
 $LG$ act in the Hilbert space $H$ through an embedding $LG \to U_{res}$ and we can simply 
 choose $T_A \equiv 1$ for all gauge connections $A.$ 
 
 In three dimensions the 
 simplest expression for the cocycle is, [Fa], [Mi85], 
 $$c(A;X,Y) =  \frac{1}{24\pi^2} \int_M \text{tr\,} A[dX, dY].\tag2.2$$

 \vskip 0.3in

 \vskip 0.3in
 
 3.  GERBAL ACTIONS AND 3-COCYCLES
 
 \vskip 0.3in
 
 Let $\Cal B$ be an associative algebra and $G$ a group. Assume that we have a group
 homomorphism $s:G \to \Out(\Cal B)$ where $\Out(\Cal B)$ is the group of outer 
 automorphims of $\Cal B,$ that is, $\Out(\Cal B) = \Aut(\Cal B)/\text{In}(\Cal B),$ all
 automorphims modulo the normal subgroup of inner automorphisms. 
 If one chooses any lift $\tilde  s: G \to \Aut(\Cal B)$ then we can write
 $$\tilde s(g) \tilde s(g') =  \sigma(g,g')\cdot \tilde s(gg')$$
  for some $\sigma(g,g')\in \text{In}(\Cal B).$ From the definition follows immediately the 
  cocycle property
  $$\sigma(g,g')\sigma(gg',g'')  = [\tilde s(g)\sigma(g', g'') \tilde s(g)^ {-1}] \sigma(g,g'g'') \text{ for all } g,g',g''\in G.\tag3.1$$ 
  
  Let next $H$ be any central extension of $\text{In}(\Cal B)$ by an abelian group $a.$ That is,
  we have an exact sequence of groups,
  $$1 \to a\to H \to \text{In}(\Cal B) \to 1.$$ 
  Let $\hat\sigma$ be a lift of the map $\sigma:G\times G  \to \text{In}(\Cal B)$ to a map $\hat\sigma: G \times G \to H.$
  We have then
  $$\hat\sigma(g,g') \hat\sigma(gg',g'') =  [\tilde s(g)\hat \sigma(g', g'') \tilde s(g)^ {-1}]   \hat \sigma(g,g'g'') 
  \cdot \alpha(g,g',g'') \text{ for all } g,g',g''\in G$$
  where $\alpha: G\times G \times G \to a .$ Here the action of the outer automorphism $s(g)$ on $\hat\sigma(*)$ is defined by
  $ s(g) \hat\sigma(*) s(g) ^{-1} =$ the lift of $ s(g)\sigma(*)
  s(g)^{-1}\in \text{In}(\Cal B)$ to an element in $H.$ One can show that $\alpha$ is a 3-cocycle [S. MacLane, Lemma 8.4],
   $$ \alpha(g_2,g_3,g_4) \alpha(g_1g_2,g_3,g_4)^{-1} \alpha(g_1,g_2g_3,g_4) \alpha(g_1,g_2,g_3g_4)^{-1}\alpha(g_1,g_2,g_3)=1.$$

 Next we construct an example from quantum field theory. 
 Let  $G$ be a compact Lie group and $P$ the space of smooth paths $f:[0,1] \to G$ with
 initial point $f(0)= e,$ the neural element, and quasiperiodicity condition $f^{-1} df$ a smooth
 function. 
 
 $P$ is a group under point-wise multiplication but it is also a principal $\Omega G$ bundle
 over $G.$ Here $\Omega G \subset P$ is the loop group with $f(0)=f(1)=e$ and $\pi: P
 \to G$ is the projection to the end point $f(1).$  Fix an unitary representation $\rho$ of $G$ in
 $\Bbb C^N$ and denote $H=L^2(S^1,\Bbb C^N).$ 
 
 For each polarization $H=H_-\oplus H_+$ we have a vacuum representation of the
 CAR algebra $\Cal B(H)$ in a Hilbert space $\Cal F(H_+).$ Denote by $\Cal C$ the category of these
 representations. Denote by $a(v), a^*(v)$ the generators 
 of $\Cal B(H)$ corresponding to a vector $v\in H,$
 $$ a^*(u) a(v) + a(v) a^*(u) = 2 <v,u>$$
 and all the other anticommutators equal to zero. 
 
 Any element $f\in P$ defines a unique  automorphism of $\Cal B(H)$
 with $\phi_f( a^*(v)) = a^*( f\cdot v),$  where $f\cdot v$ is the function on the circle defined by
 $\rho(f(x)) v(x).$ These automorphims are in general not inner except when $f$ is periodic.
 We have now a map $s: G \to \Aut(\Cal B)/\text{In}(\Cal B)$ given by $g\mapsto F(g)$ where $F(g)$ is
 an arbitrary smooth quasiperiodic function on $[0,1]$ such that $F(g)(1)=e.$ Any two such functions
 $F(g), F'(g)$ differ by an element $\sigma$ of $\Omega G,$ $F(g)(x) =F'(g)(x) \sigma(x).$ Now $\sigma$ is an inner 
 automorphism through a projective representation of the loop group $\Omega G$ in $\Cal F(H_+).$ 
 
 In an open neighborhood $U$ of the neutral element $e$ in $G$ we can fix in a smooth way
 for any $g\in U$ a 
 path $F(g)$ with $F(g)(0)=e$ and $F(g)(1) = g.$  Of course, for a connected group $G$ we can make 
 this choice globally on $G$ but then the dependence of the path $F(g)$ would not be a 
 continuous function of the end point.   For a pair $g_1, g_2\in G$ we have $\sigma(g_1,g_2)F(g_1g_2)=
 F(g_1)F(g_2)$  with $\sigma(g_1, g_1) \in \Omega G.$

 For a triple of elements $g_1, g_2,g_3$ we  have now 
 $$F(g_1)F(g_2)F(g_3) = \sigma(g_1,g_2)F(g_1g_2)F(g_3)=\sigma(g_1,g_2)\sigma(g_1g_2,g_3)F(g_1g_2g_3).$$
 In the same way,
 $$\align F(g_1)F(g_2) F(g_3)&= F(g_1)\sigma(g_2,g_3) F(g_2g_3)= [g_1 \sigma(g_2,g_3) g_1^{-1}] F(g_1)F(g_2g_3)\\
& = [g_1 \sigma(g_2,g_3) g_1^{-1}] \sigma(g_1, g_2g_3) F(g_1g_2g_3)\endalign$$
 which proves the cocycle relation (3.1).

 Lifting the loop group elements $\sigma$ to inner automorphims $\hat \sigma$ through a projective
 representation of $\Omega  G$  we can write
 $$ \hat \sigma(g_1,g_2) \hat\sigma(g_1g_2,g_3) =\text{Aut}(g_1)[\hat\sigma(g_2,g_3)]\hat\sigma(g_1,g_2g_3)
 \alpha(g_1,g_2,g_3),$$
 where $\alpha:  G\times G\times G \to S^1$ is some phase function arising from the fact that the
 projective lift is not necessarily a group homomorphism. 
 
 An equivalent point of view to the construction of  the 3-cocycle $\alpha$ is this:  We are trying to construct
 a central extension $\hat P$ of the group $P$ of paths in $G$ (with initial point $e\in G$) as an extension of the central
 extension over the subgroup $\Omega G.$ The failure of this central extension is measured by the 
 cocycle $\alpha,$ as an obstruction to associativity of  $\hat P.$ 
 On the Lie algebra level, we have a corresponding cocycle $c_3=d\alpha$ which is easily computed.
 The cocycle $c$ of $\Omega \frak{g}$ extends to the path Lie algebra $P\frak{g}$ as
 $$ c(X,Y) = \frac{1}{4\pi i} \int_{[0,2\pi]}  \text{tr\,} (XdY -YdX).$$ 
 This is an antisymmetric bilinear form on $P\frak{g}$ but it fails to be a Lie algebra 2-cocycle. The 
 coboundary is given by 
 $$ \align (\delta c)(X,Y,Z) & =c(X,[Y,Z]) + c(Y,[Z,X]) + c(Z,[X,Y]) \\
 &= -\frac{1}{4\pi i} \text{tr\,} X[Y,Z]\vert_{2\pi} =d\alpha(X,Y,Z).\endalign$$
 Thus $\delta c$ reduces to a 3-cocycle of the Lie algebra $\frak{g}$ of $G$ on the boundary 
 $t=2\pi.$  This cocycle defines by left translations on $G$ the left-invariant de Rham form
 $-\frac{1}{12\pi i}\text{tr\,}(g^{-1} dg)^3;$ this is normalized as $2\pi i$ times an integral 
 3-form on $G.$ 
 
 Let $f_1,f_2\in P$ and $f_{12}\in P$ with the property $f_1(2\pi)f_2(2\pi) = f_{12}(2\pi).$ Then we have a fiber $S^1$ over the loop 
 $\phi_{12}= f_1(t) f_2(t) f_{12}(t)^{-1}$ coming from the central extension $\widehat{\Omega G} \to \Omega G.$ Concretely, this fiber can be realized 
 geometrically as a pair $(f,\lambda)$ where $f:D \to G,$   $D$ is the unit disk with boundary $S^1$ such that the restriction of $f$ to $S^1$
 is the loop above,  and $\lambda\in S^1.$ Two pairs $(f,\lambda), (f',\lambda')$ are equivalent if the restrictions of $f,f'$ to the boundary
 are equal and $\lambda' = \lambda e^{W(f,f')}$ where 
 $$W(f,f') = \frac{1}{12\pi i} \int_B \text{tr\,} (g^{-1} dg)^3,$$
 where $B$ is a unit ball with boundary $S^2= S^2_+ \cup S^2_-$ and $g$ is any extension to $B$ of the map $f\cup f'$ on $S^2$ 
 obtained by joining $f,f'$ on the boundary circle $S^1$ of the unit
 disk $D.$ The product in the central extension of the full loop group
 $LG$ is then defined as 
 $$ [(f,\lambda)] \cdot [(f',\lambda')] = [(ff', \lambda\lambda' e^{\gamma(f,f')})], $$ 
 where
 $$ \gamma(f,f')= \frac{1}{4\pi i} \int_D \text{tr\,} f^{-1}df \wedge df' {f'}^{-1},$$
 see [Mi87] for details; here the square brackets mean equivalence
 classes of pairs, subject to the equivalence defined above. 
 
 The 3-cocycle $\alpha$ can now be written in terms of the local data as
 $$\alpha = \exp[  \gamma(\phi_{12}, \phi_{12,3})- \gamma(\text{Aut}_{f_1} \phi_{23},\phi_{1,23}) +W(h) ]$$
 where $\phi_{12,3}$ is an extension to the disk of the loop composed from the paths $f_1f_2$ and $f_3$
 and the path $f_{12,3}$ joining the identity $e$ to $g_1g_2g_3.$  $h$ is the function on $D,$ equal to the neutral element
 on the boundary, such that $\phi_{12} \phi_{12,33} = \phi_{12} \phi_{12,3}h.$ The value of $\alpha$ depends now,
 besides on the paths $f_i,$ on the extensions $\phi$ to the disk $D$ of the boundary loops determined
 by $f_i, f_{12}, f_{23},f_{12,3}, f_{1,23}.$ However, different choices of extensions are related by phase factors
 which can be obtained from the equivalence relation
 $$(\phi,\lambda) \equiv (\phi h, e^{\gamma(\phi,h)+ W(h) })$$ 
 defining the central extension of the loop group.

 \vskip 0.3in
 
 4. THE CASE OF HIGHER DIMENSIONS
 
 \vskip 0.3in
 
 The construction of the gerbal action of $G$ has a generalization which comes from a study of gauge 
 anomalies in higher dimensions. Fix again a compact Lie group  $G$ and denote now by $ MG$ the group of
 smooth maps from a compact manifold $M$ to $G,$ which is an infinite-dimensional Lie group under point-wise 
 multiplication of maps.  Assume also that $M$ is 
 a boundary of a compact manifold $N$ and denote by $NG$ the group of maps from $N$ to $G$ such that 
 the normal derivatives at the boundary $M$ vanish in all orders. Finally, let $\Cal G$ be the normal  subgroup of $NG$ 
 consisting of maps equal to the constant $e\in G$ on the boundary. this will play the role of $\Omega G$ in the previous section.
 Now $NG/\Cal G = MG.$ 
 
 We also assume that $M$ is a boundary of a contractible manifold $N'$ and denote by $\overline N$ the manifold obtained
 from $N, N'$ by gluing along the common boundary. We also assume a spin structure and Riemann metric given on $\overline N.$
 
 We may view elements $f\in NG$ as $\frak g$ valued vector potentials on the space $\overline  N.$ 
 This correspondence is given by $f\mapsto A = f^{-1} df$ on
 $N$ and $A'$ is fixed from the boundary values $A\vert_M$ and from a contraction of $N'$ to one point.
 This construction gives a map from the  group $MG$  to the moduli space of gauge connections in a trivial vector bundle over
 $\overline N.$ 
 
 \bf Example \rm Let $M=S^n,$ viewed as the equator in the sphere $S^{n+1}.$  Fix a path $\alpha$ from the South Pole of $S^{n+1}$ to the 
 North Pole. Then for any great circle joining the North Pole to the South Pole we can take the union with $\alpha$, giving a loop 
 starting from the South Pole and traveling via the North Pole. For given vector potential $A$ let $g_A$ be the holonomy around this loop.
 The great circles are parametrized by points on the equator $S^n$ and thus we obtain a map $S^n \to g_A.$  The group $\Cal G$ of \it based \rm 
 gauge transformations, those which are equal to the identity on the South Pole, does not affect the holonomy, thus we obtain 
 a map from the moduli space $\Cal A/ \Cal G$ of gauge connections to $S^n G.$ The same construction can be made for any $n$-sphere
 $S^n$ around the South Pole, and we may view $S^{n+1}$ as the $(n+1)$-dimensional solid ball $N$ with boundary $M=S^n$ contracted to
 one point, the South Pole of $S^{n+1}.$ Here $\Cal G$ is viewed as the group $G$ valued maps on $N$ equal to the constant $e$ on
 the boundary. In this case one can actually show that the gauge moduli space $\Cal A/\Cal G$ is homotopy equivalent to $MG.$ 
 
 The main difference as compared to the case of loop group is that the transformations $f\in \Cal G$ are not in general implementable 
 Bogoliubov automorphisms. However, as explained in Section 2,   they are well-defined automorphism of 
 a Hilbert bundle $\Cal F$ over $\Cal A,$ the space of $\frak{g}$ valued connections on $\overline N.$ 
 The fibers of this bundle are fermionic Fock spaces, each of them carries an 
 (inequivalent) representation of the canonical anticommutation relations, with a Dirac vacuum which depends on the background field
 $A\in \Cal A.$ The group $\Cal G$ acts on this bundle through an abelian extension $\widehat{\Cal G}.$ 
 
 The 3-cocycle is constructed in essentially the same way as in Section 3. So for an element $g\in MG$ select $f\in NG$ such that 
 the restriction of $f$ to the boundary is equal to $g.$ For any pair $g_1,g_2\in MG$ we have then $f_1f_2 f_{12}^{-1}\in \Cal G$ where 
 again $f_{12} \in NG$ such that $f_{12} \vert_M = g_1g_2.$ For a triple $g_1,g_2,g_3 \in MG$ we then construct 
 $\alpha(A; g_1,g_3,g_3) \in S^1$ as before, but now it depends on the connection $A$ since the operators $\hat\sigma$
 now all depend on $A.$ 
 
  $$\hat \sigma_{A'}(g_1,g_2) \hat\sigma_A(g_1g_2,g_3) =\text{Aut}(g_1)[\hat\sigma_{A''}(g_2,g_3)]\hat\sigma_A(g_1,g_2g_3)
 \alpha(A;g_1,g_2,g_3),$$
 where $A'= A^{\sigma(g_1g_2,g_3)}$ is the gauge transform of $A$ by $\sigma(g_1g_2,g_3) \in\Cal G$ and similarly
 $A''$ is the gauge transform of $A$ by $\sigma(g_1,g_2g_3).$

 \bf Example \rm  Again, passing to the Lie algebra cocycles one gets reasonably simple expressions. For example, in the
 case dim$M = 2$  the Lie algebra extension of Lie$(\Cal G)$ is given by the 2-cocycle (2.2) and for a manifold $N$ with boundary 
 $M$ this formula is not a cocycle but its coboundary is the Lie algebra 3-cocycle 
 $$ d\alpha(X,Y,Z) = - \frac{1}{8\pi^2} \int_M \text{tr\,} X[dY,dZ] . \tag4.1$$ 
 In this case the cocycle does not depend on the variable $A$ but when dim$M >2$ it does. 
 
 \vskip 0.3in
 
  5.  TWISTED K-THEORY ON MODULI SPACES OF GAUGE CONNECTIONS
 
 \vskip 0.3in
 
 Let $P$ be a principal bundle over a space $X$ with model fiber equal to the projective 
 unitary group $PU(H)= U(H)/S^1$ of a complex Hilbert space. It is known that equivalence 
 classes of such bundles are classified by elements in $H^3(X,\Bbb Z),$ the Dixmier-Douady
 class of the bundle, [DD]. 
 
 K-theory of $X$ twisted by $P,$  denoted by $K^*(X,P),$ is defined as the abelian group of homotopy classes of
 sections of a bundle $Q,$ defined as an associated bundle with fiber equal to the ($Z_2$ graded) space of Fredholm
 operators in $H$ with $PU(H)$ action given by the conjugation $T\mapsto gTg^{-1}.$ The grading is as in ordinary complex 
 $K$ theory: The even sector is defined by the space of \it all \rm Fredholm operators whereas the odd sector is
 defined by self-adjoint operators with both positive and negative essential spectrum.  As a model, one can use
 either bounded Fredholm operators, or unbounded operators for example with the graph topology, [AtSe].
 
 If $X=G$ is a compact Lie group  one can construct elements of $K^*(G,P)$ in terms of highest weight 
 representation of the central extension $\widehat{LG},$  [Mi04 ].  Actually, these come as $G$ equivariant classes,
 under the conjugation action of $G$ on itself.  In the equivariant case the construction of $K^*(G, P)$ is related to
 the Verlinde algebra in conformal field theory, [FHT].  Although for simple compact Lie groups there exists classification 
 theorems [Do], [Br] also in the nonequivariant case it is still an open problem how to give explicit constructions for all
 classes in
 the nonequivariant case,
 in terms of families of Fredholm operators,
 using representation theory even for unitary groups $SU(n)$ when $n > 3.$

 Let $\omega: \Cal A \times \Cal G \to U_{res}(H_+\oplus H_-)$ be the 1-cocycle 
 constructed in Section 2. Let $Y$ be a family of Fredholm operators in $\Cal F$ which is
 is mapped onto itself under a projective representation $g\mapsto \hat g$ of $U_{res}$ in $\Cal F,$ 
 $T\mapsto \hat g T \hat g^{-1}\in Y$ for any $T\in Y.$ 
 
 Now we have an action of a central extension of the groupoid $(\Cal A, \Cal G)$ on $Y$ by
 $$(A,g): Y\to Y,  T\mapsto\hat\omega(A;g) T \hat\omega(A;g)^{-1}.$$
 We can also view this as a central extension of the transformation groupoid defined by the action of the 
 gauge group $\Cal G$ on the space $\Cal A \times Y.$ If $\Cal  G$ is the group of based gauge transformations
 then it acts freely on $\Cal A$ and therefore also freely on $\Cal A \times Y$. If furthermore $Y$ is 
 contractible then $(\Cal A \times Y)/ \Cal G \simeq \Cal A/\Cal  G$ is the gauge moduli space.
 
 A system of Fredholm operators transforming covariantly under $U_{res}$ can be constructed 
 from a Dirac operator on the infinite-dimensional Grassmann manifold $Gr_{res}=U_{res}(H_+\oplus H_-)/(U(H_+) \times U(H_-) ),$
 [T\"a].  The members of the family are parametrized by a gauge connection on a complex line bundle $L$ over $Gr_{res}.$
 The line bundle $L$ is used as twisting of the spin bundle over $Gr_{res},$ and can be viewed as defining 
 a spin$_{C}$ structure on $Gr_{res}.$

 In the case when $\Cal A$ is the space of gauge connections on the circle and $\Cal G= LG$ is
 a loop group the cocycle $\omega(A;g)$ does not depend on $A$ and it gives a unitary
 representation of $LG$ in the Hilbert space $H$ and $g\mapsto \hat g$ is given by a 
 representation of a central extension $\widehat{LG}$ in $\Cal F.$ As $Y$ we can take the
 family $Q_A$ of supercharges in [Mi04]  parametrized by points in $\Cal A.$ This means that in the notation above,
 we can identify $\Cal A$ as the diagonal in $\Cal A \times Y$ and we have a natural identification of the 
 groupoid moduli space with the moduli space of gauge connections on the circle.
 Here the groupoid action
 defines an element in the twisted $G$-equivariant K-theory  on $G.$   In fact, in the case of the circle, one
 can directly work with highest weight representations of the loop group without using the embedding
 $LG \subset U_{res}.$

\vskip 0.3in

REFERENCES

\vskip 0.3in

 [ABJ]  S.L. Adler: Axial-vector vertex in spinor electrodynamics. Phys. Rev, \bf 177, \rm p. 2426 (1969).
 W.A. Bardeen: Anomalous Ward identities in spinor field theories. Phys. Rev. \bf 184, \rm p.1848 (1969)
   J. Bell and R. Jackiw: The PCAC puzzle: $\pi^0 \to \gamma\gamma$ in the $\sigma$-model. Nuevo Cimento  \bf A60, \rm 
   p. 47 (1969). 
 
 [AS]  M.F. Atiyah and I.M. Singer: Dirac operators coupled to vector potentials. Proc. Natl. Acad. Sci. USA,
 \bf 81, \rm pp. 2597-2600 (1984)
 
 [AtSe]  M. Atiyah and G. Segal: Twisted $K$-theory.  Ukr. Mat. Visn.
 \bf 1, \rm ,  no. 3, 287--330 (2004);  translation in
 Ukr. Math. Bull. 
 \bf 1, \rm   no. 3, pp. 291--334 (2004)
 
 [Br] V. Braun:  Twisted $K$-theory of Lie groups.  J. High Energy Phys.   no. 3, 029, 15 pp (2004) (electronic)
 
 [Ca]  A. L. Carey: The origin of three-cocycles in quantum field theory. Phys. Lett. B \bf 194, \rm pp. 267-270 (1987)
 
 [CGRS] A.L. Carey, H. Grundling, I. Raeburn, and C. Sutherland: Group
 actions on $C^*$-algebras, 3-cocycles and quantum field
 theory. Commun. Math. Phys. \bf 168, \rm pp. 389-416 (1995)
 
 [DD] J. Dixmier and A. Douady:  Champs continus d'espaces hilbertiens et de $C\sp{\ast} $-algbres.   Bull. Soc. Math. France  \bf 91, \rm  pp. 227--284 (1963)
 
 [Do] C. Douglas: On the twisted $K$-homology of simple Lie groups.  Topology  \bf 45, \rm   pp. 955--988 (2006)
 
 [Fa] L.D. Faddeev: Operator anomaly for the Gauss law. Phys. Lett. B \bf 145, \rm p.81 (1984)
 
 [FHT]  D. Freed, M. Hopkins, and C. Teleman: Twisted K-theory and loop group representations, arXiv:math/0312155. 
 Twisted equivariant $K$-theory with complex coefficients.  J. Topol.  1  (2008),  no. 1, 16--44. 
 
 [FZ]  E. Frenkel and Xinwen Zhu: Gerbal representations of double loop groups. \newline
  arXiv.math/0810.1487
 
 [GJJ] B. Grossman. The meaning of the third cocycle in the group cohomology of nonabelian gauge theories. 
 Phys. Lett. B  \bf 160, \rm   pp.  94--100 (1985).
 R. Jackiw, R: Three-cocycle in mathematics and physics.  Phys. Rev. Lett.  \bf 54, \rm pp.159--162 (1985).
 S.G. Jo: Commutator of gauge generators in nonabelian chiral theory.  Nuclear Phys. B  \bf 259, \rm   pp. 616--636 (1985)
 
 [LM]  E. Langmann and  J. Mickelsson:  Scattering matrix in external field problems.  J. Math. Phys.  \bf 37, \rm  pp. 3933--3953 (1996)
 
 [Mi85] J. Mickelsson:  Chiral anomalies in even and odd dimensions.  Comm. Math. Phys.  \bf 97  (1985),  \rm pp. 361--370 (1985)
 
 [Mi87] J. Mickelsson: Kac-Moody groups, topology of the Dirac determinant bundle, and fermionization.  Comm. Math. Phys.  \bf 110, \rm pp. 173--183 (1987)
 
 [Mi90] J. Mickelsson:  Commutator anomalies and the Fock bundle.  Comm. Math. Phys.  \bf 127, \rm  pp. 285--294
 (1990)
 
 [Mi93] J. Mickelsson:  Hilbert space cocycles as representations of $(3+1)$-D current algebras.  Lett. Math. Phys.  \bf 28, \rm  pp.97--106 (1993).
 Wodzicki residue and anomalies of current algebras.  Integrable models and strings (Espoo, 1993),  pp. 123--135, Lecture Notes in Phys., \bf 436, \rm
 Springer, Berlin (1994)
 
 [Mi04]  J. Mickelsson: Gerbes, (twisted) $K$-theory, and the supersymmetric WZW model.  Infinite dimensional groups and manifolds,  93--107, IRMA 
 Lect. Math. Theor. Phys., \bf 5, \rm de Gruyter, Berlin (2004)
 
 [ML] Saunders Mac Lane: \it Homology. \rm Die Grundlehren der Mathematischen Wissenschaften, Band 114.
    Springer Verlag (1963)
    
 [PS] A. Pressley and  G. Segal: \it Loop Groups. \rm Clarendon Press, Oxford (1986)
 
 [T\"a] V. T\"ahtinen: The Dirac operator on restricted Grassmannian, in progress

 \enddocument